\documentclass[10pt]{amsart}
\usepackage{amssymb,amsthm,amsmath,amstext}
\usepackage{stmaryrd}
\usepackage{mathrsfs} 
\usepackage{float}

\usepackage{graphicx}
\usepackage{amsfonts}

\usepackage[all,cmtip]{xy}
\usepackage{amssymb}
\usepackage{enumitem}
\usepackage{booktabs}
\usepackage{longtable}
\usepackage{caption}

\newtheorem{theorem}{\bf Theorem}[section]
\newtheorem{definition}[theorem]{\bf Definition}
\newtheorem{corollary}[theorem]{\bf Corollary}
\newtheorem{lemma}[theorem]{\bf Lemma}
\newtheorem{proposition}[theorem]{\bf Proposition}
\newtheorem{example}[theorem]{\bf Example}
\newtheorem{remark}[theorem]{\bf Remark}

\newtheorem*{theorem*}{Theorem}

\usepackage[left=1.5in,top=1.2in,right=1.5in,bottom=1.2in]{geometry}

\usepackage[backref=page]{hyperref}
\usepackage{hyperref}
\hypersetup{pdftitle={Submanifold with splitting tangent sequence}}
\hypersetup{pdfauthor={Cong Ding}}
\hypersetup{colorlinks=true,linkcolor=blue,anchorcolor=blue,citecolor=blue}

\begin{document}
	
	\title[Splitting submanifold]{Splitting submanifolds in rational homogeneous spaces of Picard number one
	}
	
	\author{Cong Ding}
	\address{Department of Mathematics, The University of Hong Kong \\Pokfulam Road, Hong Kong, China}
	\address[Current address]{Morningside Center of Mathematics, Academy of Mathematics \& Systems Science, The Chinese Academy of Sciences, Beijing, China}
	\email{congding@amss.ac.cn}
	
	
\maketitle

\begin{abstract}
	Let $M$ be a complex manifold. We prove that a compact submanifold $S\subset M$ with splitting tangent sequence (called a splitting submanifold) is rational homogeneous when $M$ is in a large class of rational homogeneous spaces of Picard number one. Moreover, when $M$ is 
irreducible Hermitian symmetric, we prove that $S$ must be also Hermitian symmetric. These cover some of the results given in \cite{MR2190340}. The basic tool we use is the restriction and projection map $\pi$ of the global holomorphic vector fields on the ambient space which is induced from the splitting condition. The usage of global holomorphic vector fields may help us set up a new scheme to classify the splitting submanifolds in explicit examples, as an example we give a differential geometric proof for the classification of compact splitting submanifolds with $\dim\geq 2$ in a hyperquadric, which has been proven in \cite{MR2190340} using algebraic geometry.

\end{abstract}

	\section{Introduction}
\label{intro}
Let $M$ be a complex manifold and $S$ be a complex submanifold, there is a tangent sequence \[0\longrightarrow T(S) \longrightarrow T(M)|_S \longrightarrow N_{S|M}\longrightarrow 0\]
If this short exact sequence splits holomorphically over $S$, i.e. $T(M)|_S=T(S)\oplus N_{S|M}$, we will say that $S$ is a submanifold with splitting tangent sequence and for short we call it \textbf{a splitting submanifold}. Note that the splitting condition is purely a holomorphic condition and a prior not related to the K\"{a}hler metric. And in fact it is weaker than total geodesy in general, as total geodesy is equivalent to orthogonality of the direct sum with respect to the K\"{a}hler metric. For example, we consider a hyperquadric $Q^n$ and a geodesic 2-sphere $\mathbb{P}^1\times \mathbb{P}^1 \subset Q^n$, then the graph of a holomorphic map of degree $2$ from $\mathbb{P}^1$ to $ \mathbb{P}^1$ gives a splitting curve in $Q^n$ which is not geodesic with respect to any K\"ahler-Einstein metric on $Q^n$. Also we can observe that splitting condition does not have direct relation with equivariance, for example, the natural equivariant embedding $Q^n \subset \mathbb{P}^{n+1} \subset Q^m$ is not splitting since $Q^n$ is not splitting in $\mathbb{P}^{n+1}$, on the other hand the graph example for splitting curves tells that splitting condition does not imply equivariance.

\begin{remark}
	Although splitting condition is weaker than total geodesy, we have a characterization of splitting condition using Hermitian metric instead of K\"{a}hler metric. More precisely we have $S\subset M$ is splitting if and only if there is a Hermitian metric on $M$ such that there is an orthogonal decomposition $T(M)|_S=T(S)\oplus N_{S|M}$ with respect to the Hermitian metric. One direction is straightforward, for another direction if $S$ is splitting in $M$ and $T(M)|_S=T(S)\oplus N_{S|M}$ is a direct sum of $T(S)$ and its complementary bundle, we can choose a Hermitan metric $h$ on $T(M)|_S$ such that the direct sum is orthogonal with respect to $h$. Then metric can be extended to an open neighborhood of $S$ and by the argument of partition of unity it can be extended to a Hermitian metric on $M$.\end{remark}

General purpose of the research is to study the properties or even classifications of splitting submanifolds in certain Hermitian (or locally Hermitian) symmetric spaces and more generally rational homogeneous manifolds. Within this article, the submanifold $S$ is assumed to be compact if not specified. Van de Ven\cite{MR0116361} obtained that the splitting submanifolds in $\mathbb{P}^n$ are linear thus totally geodesic with respect to any Fubini-Study metric. A generalization to all K\"{a}hler manifolds of constant holomorphic sectional curvature can be found in Mok\cite{MR2156495} where harmonic forms of the obstruction class of tangent sequence were used.

Jahnke\cite{MR2190340} considered more general cases when the ambient space $M$ is some irreducible Hermitian symmetric space of compact type and obtained that

\begin{theorem}\label{thm1}\rm(\textbf{Jahnke, \cite[Theorem~4.7]{MR2190340}})
	Any splitting submanifold with $\dim\geq2$ in a hyperquadric $Q^{n}$ is a linear subspace or a smooth complete linear section subquadric. Any splitting curve in $Q^{n}$ is rational.
\end{theorem}
\begin{theorem}\label{thm2}\rm(\textbf{Jahnke, \cite[Theorem~5.2]{MR2190340}})
	Any splitting submanifold in a Grassmannian is rational homogeneous.
\end{theorem}

The splitting condition will quickly tell that the submanifold $S$ is homogeneous since its tangent bundle is globally generated. From the structure theorem of Borel and Remmert for compact homogeneous K\"{a}hler manifolds (cf.\cite{MR0145557}), we have a decomposition $S=A\times N$ where $N$ is a rational homogeneous space and $A$ is an abelian variety. The first task is to rule out the abelian variety part and from Theorem~\ref{thm1} and Theorem \ref{thm2} we know that it does not exist when the ambient space is a hyperquadric or a Grassmannian. In general, Jahnke provided an upper bound for the dimension of possible abelian variety related to some notions of positivity of tangent bundles, i.e. $k$-ampleness and ampleness in the usual sense (cf.\cite{MR2190340}).

\begin{proposition}\label{prop1}\rm(\textbf{Jahnke, \cite[Proposition~2.2]{MR2190340}})
	Let $S$ be a splitting submanifold in a rational homogeneous space $M$. Let $k,l$ be the minimal numbers such that $T(M)$ is $k$-ample and $\bigwedge^{l}T_M$ is ample. If dim$(S)\geq$ min$(k+1,l)$, then $S$ is rational homogeneous.
\end{proposition}

In this article we will prove that the abelian variety part does not exist for lots of rational homogeneous spaces of Picard number one including the cases when $M$ is irreducible Hermitian symmetric, so it covers the Grassmannian and hyperquadric case.

Moreover, when the ambient space $M$ is an irreducible Hermitian symmetic space of compact type, it is natural to ask whether $S$ is Hermitian symmetric. This is motivated by the work of Mok-Ng\cite{MR3647134} where the answer is affirmative for the quotient of its noncompact dual when $\dim(S)$ is big enough, and the lower bound of $\dim(S)$ is related to the curvature tensor in their results.

For compact type, to determine whether a rational homogeneous space can be splitting in an irreducible compact Hermitian symmetric space is not obvious, even for the simple example Symplectic Grassmannian embedded in a Grassmannian. In this article we will prove such $S$ must be biholomorphic to a compact Hermitian symmetric space, so no non-symmetric rational homogeneous spaces can be splitting in irreducible compact Hermitian symmetric spaces, this seems more rigid than noncompact type case, as there is no dimension restriction for the submanifold $S$.

We know the splitting condition is a global condition. Our method throughout the paper is based on the restriction and projection map $\pi$ of global holomorphic vector fields of the ambient space, i.e. for every point $x\in S$ and every global holomorphic vector field $Z$ on $M$, $\pi(Z)_{x}$ is the projection from $T_x(M)$ to $T_{x}(S)$, induced by the transversality $T_x(M)=T_{x}(S)\oplus N_{S|M,x}$.

To prove the rational homogeneity, we will use global holomorphic vector fields of $M$ vanishing at the reference point $o$. Moreover, to prove the Hermitian symmetry, we will use the Euler vector field, which is a global holomorphic vector field vanishing at the reference point to the first order. Motivated by this we furtherly consider the 
restriction and projection of global holomorphic vector fields vanishing at the reference point to the second order. From this we can obtain some restrictions on the tangent subspaces and its complement pointwisely. Therefore we may obtain a possible procedure to give the full classification of splitting submanifolds for explicit examples: (i) Classify the tangent subspaces (Infinitesimal information). (ii) Prove the integrability (Global information). As an example, we give a new proof on the classification of compact splitting submanifolds with $\dim\geq 2$ in a hyperquadric (see Theorem \ref{thm1}). As an attempt on classification for splitting submanifolds 
in Grassmannians, we follow this procedure to study splitting submanifolds with $\dim\geq 2$ in the low-dimensional case $M=G(2,3)$ in \cite{DingThesis}.

\section{Statement of Main results}
\label{MainResults}
On rational homogeneity we can prove

\begin{theorem}\label{rational splitting}
	Let $M=G/P$ be a rational homogeneous space of Picard number one, where $G$ is a connected complex simple Lie group and $P$ is a maximal parabolic subgroup. Suppose $S\subset M$ is a compact submanifold with splitting tangent sequence. Then $S$ is rational homogeneous if depth of the gradation for $\mathfrak{g}=Lie(G)$ is not greater than $2$.
\end{theorem}

\begin{remark}
	For Theorem \ref{rational splitting}, it covers a large class of rational homogeneous spaces including all classical types. In fact similar method may apply to other cases not included above with higher depth, a computation for $(E_6,\alpha_4)$ (depth $=3$) can be found in \cite{DingThesis}. However it requires a tedious combinatorial checking on root systems, other cases seems to be too lengthy to finish, a more elegant method may be required. 
\end{remark}

Furthermore when $M$ is symmetric we have 
\begin{theorem} 
	If $S$ is a compact splitting submanifold in an irreducible compact Hermitian symmetric space $M$, then $S$ is biholomorphic to a Hermitian symmetric space.
\end{theorem}

\section{Preliminaries}
\label{Prel}

\subsection{Rational homogeneous spaces}
\label{RHS}

We give some basic introduction about rational homogeneous spaces, which is a large class of homogeneous spaces including compact Hermitian symmetric spaces. For general reference, we refer \cite{MR1274961} or section 2 of \cite{MR1681093}.

\begin{definition}
	A complex homogeneous space is called a rational homogeneous space if it is projective algebraic.
\end{definition}

Any rational homogeneous space can be written as a product of irreducible rational homogeneous spaces. When $M$ is irreducible we can write $M=G/P$ where $G$ is a connected complex simple Lie group and $P$ is a parabolic subgroup. $M$ is of Picard number one if and only if $P$ is maximal parabolic. Let $\mathfrak{g}$ be the Lie algebra of $G$ and $\mathfrak{p}$ be the Lie algebra of $P$. Fix a Levi decomposition $\mathfrak{p}=\mathfrak{l}^{\mathbb{C}}+\mathfrak{u}$ where $\mathfrak{l}^{\mathbb{C}}$ is reductive and $\mathfrak{u}$ is nilpotent. Choose a Cartan subalgebra $\mathfrak{t}^{\mathbb{C}}\subset \mathfrak{l}^{\mathbb{C}}$ and a root system $\Delta \subset (\mathfrak{t}^{\mathbb{C}})^{*}$, we use $\mathfrak{g}_{\alpha}$ to denote the root space corresponding to $\alpha \in \Delta$. Fix a system of simple roots $\Psi=\{\alpha_1,...\alpha_{n}\}$ and choose a set of distinguished simple roots $\Psi'=\{\alpha_{i_1},..., \alpha_{i_\ell} \}$ where $\ell$ is the Picard number of $M$, then define $\Delta_{k}=\{\alpha \in \Delta, \alpha=\sum^{n}_{q=1}m_{q}\alpha_{q},\sum_{r=i_1}^{i_\ell}m_{r}=k\}$ where $-m\leq k\leq m$ and $m$ is the largest number with $\Delta_{m}\neq 0$. This gives a graded structure on $\mathfrak{g}$ such that $\mathfrak{g}=\oplus_{k=-m}^{m}\mathfrak{g}_{k}$ where $\mathfrak{g}_{0}=\mathfrak{t}^{\mathbb{C}}\oplus\bigoplus_{\alpha \in \Delta_{0}}\mathfrak{g}_{\alpha}$ and $\mathfrak{g}_{k}=\bigoplus_{\alpha \in \Delta_{k}}\mathfrak{g}_{\alpha}$, $m$ is also called the depth of the gradation. Then $\mathfrak{l}^{\mathbb{C}}=\mathfrak{g}_0$ and $\mathfrak{p}=\oplus_{k=-m}^{0}\mathfrak{g}_{k}$. For the gradation there exists a central element $\eta_0$ in $\mathfrak{g}_0$ such that $[\eta_0,\xi_k]=k\xi_k$ for any $\xi_k\in \mathfrak{g}_k$. In this case we say $M$ is \textbf{a rational homogeneous space associated to $\Psi'$}. When $\Psi'=\{\alpha_i\}$, i.e. $M$ is of Picard number one, we can use $(\mathfrak{g}, \alpha_i)$ to denote $M$. Also, we identify the holomorphic tangent space of $M$ at the reference point $o$ with the positive part in the gradation, i.e. $T_o(M) \cong \oplus_{k=1}^{m}\mathfrak{g}_{k} $. Especially, when $m=1$, then $\ell=1$ and $M$ is Hermitian symmetric and when $m=2$ and dim$\mathfrak{g}_{2}=1$, $M$ is contact homogeneous. A typical example of non-symmetric rational homogeneous space is the Symplectic Grassmannian, i.e. all the isotropic $p$-planes in a $2n$-dimensional complex symplectic vector space, where $p<n$.

\subsection{Hermitian symmetric spaces}
\label{HSS}

Inheriting the notations above, we know when $m=1$, $M$ is irreducible Hermitian symmetric and in this case , $M$ can be written as $G/P=G_c/K$ where $G_c$ is a compact real form of $G$ and $K$ is a maximal compact subgroup in $G_c$. Let $\mathfrak{l}=Lie(K)$,  and the complexification is denoted by $\mathfrak{l}^{\mathbb{C}}$. In this case the graded structure is identified with the Harish-Chandra decomposition $\mathfrak{g}=\mathfrak{l}^{\mathbb{C}}\oplus\mathfrak{m}^+\oplus \mathfrak{m}^-$ where $\mathfrak{l}^{\mathbb{C}}=\mathfrak{g}_0$, $\mathfrak{m}^+=\mathfrak{g}_1\cong T_o(M)$ and $\mathfrak{m}^-=\mathfrak{g}_{-1}$. We use $\Delta_{M}^+(\Delta_{M}^-)$ to denote the set of roots with root space contained in $\mathfrak{g}_1(\mathfrak{g}_{-1})$, which are called noncompact positive (negative) roots. Let $\mathfrak{g}_c=Lie(G_c)$. Since the Killing form on $\mathfrak{g}_c$ is negative definite, it induces an invariant Hermitian metric on $M$ which is in fact K\"{a}hler. Now the curvature tensor is denoted by $R$ and the holomorphic bisectional curvature is nonnegative.
\par In Hermitian symmetric case we can construct a maximal strongly orthogonal roots $\Pi\subset\Delta^+_{M}$ starting from the highest root. Write $\Pi=\{\beta_1,...,\beta_r\}$, then $r$ is the rank of $M$. Now we introduce the Polydisk and  Polysphere Theorem, 

\begin{theorem}[The Polydisk and Polysphere Theorem, cf.\cite{MR0404716} or \cite{MR1081948}]
	Suppose that $(X_0,g_0)$ is a rank $r$ irreducible Hermitian symmetric space of noncompact type, $X_0=G_0/K$. Then there exists a totally geodesic complex submanifold $(D,g_0|_D)$ isometric to a product of $r$ disks equipped with Poincar\'{e} metric and the $K$-action on $D$ exhausts $X_0$, i.e. $X_0=\bigcup_{k\in K}kD$. And for the dual case $(X_c, g_c)$ there exists a totally geodesic submanifold $S$ in $X_c$ isometric to a product of $r$ Riemann spheres equipped with Fubini-Study metric. Moreover $D$ is contained in $S$ through the Borel embedding of $X_0$ in $X_c$. 
\end{theorem}

We set up the following normalization for root vectors. For each root $\varphi$, we denote the coroot by $h_\varphi$ and we choose a normalized vector $e_\varphi, e_{-\varphi}$ such that $[e_\varphi,e_{-\varphi}]=H_{\varphi}$ where $H_{\varphi}$ is the dual element of $\varphi$ in $\mathfrak{t}^{\mathbb{C}}$ with respect to the Killing form of $\mathfrak{g}$ (denoted by $B$) and $h_{\varphi}=2\frac{H_\varphi}{B(H_\varphi,H_\varphi)}$, $\mathbb{C}\{e_\varphi, e_{-\varphi}, H_{\varphi}\}$ gives a three dimensional simple subalgebra isomorphic to $\mathfrak{sl}(2,\mathbb{C})$, in this normalization $B(e_\varphi,e_{-\varphi})=1$.
The holomorphic tangent space of the polysphere 
at $o$ is spanned by $\{e_{\beta_1},...,e_{\beta_r}\}$. Together with the Restricted Root Theorem (cf.\cite{MR0161943}) we have the following classification for the normal forms (representatives) for the $K$-orbits and $P$-orbits on the projectivized tangent space of $M$. 
\begin{proposition}\label{HSS-normal}
	For a nonzero holomorphic tangent vector $v\in T_o(M)$, there exists some $k\in K$ such that $Adk.v=\sum_{i=1}^{r_0}a_ie_{\beta_i}$ with $a_1\geq a_2 \cdots \geq a_{r_0}>0$ for some $1\leq r_0\leq r$. Moreover by $\mathbb{C}^*$-action induced by the Cartan subalgebra $\mathfrak{t}^{\mathbb{C}}$,
	there exists some $p\in P$ such that $Adp.v=\sum_{i=1}^{r_0}e_{\beta_i}$. Here $r_0$ is called the rank of $v$.
\end{proposition}

This provides some motivation for us to classify normal forms for the $P$-orbits on the projectivized tangent space of a general rational homogeneous space of Picard number one, which will be useful in proving the rational homogeneity of the splitting submanifolds. 

For our own use, we discuss some basic VMRT theory for irreducible compact Hermitian symmetric space $M$ without introducing the general theory developed by Hwang-Mok (see for example \cite{MR1748609},\cite{MR2521656}). Since the Picard group $Pic(M)\cong \mathbb{Z}$, we can use $\mathcal{O}(1)$ to denote the positive generator for $Pic(M)$, and the space of global holomorphic sections $\Gamma(M,\mathcal{O}(1))$ gives the so-called minimal embedding of $M$ into $\mathbb{P}(\Gamma(M,\mathcal{O}(1))^*)$. Minimal rational curves are free rational curves of minimal degree with respect to $\mathcal{O}(1)$ which are the projective lines through the minimal embedding. Fix a point $o\in M$, the variety of minimal rational tangents (VMRTs for short) denoted by $\mathcal{C}_o(M) \subset \mathbb{P}T_o(M)$ is the collection of all tangent vectors of free minimal rational curves passing through $o$ (which are also called characteristic tangent vectors). For example we know for hyperquadrics $\mathcal{C}_o(Q^n)=Q^{n-2}$. A tangent subspace is called characteristic if all tangent vectors in it are characteristic. For irreducible Hermitian symmetric spaces of compact type, there are equivalent characterization for minimal rational tangents, they are highest weight vectors with respect to the isotropy action and also the tangent vectors with maximal holomorphic sectional curvature (normalized to unit length with respect to the K\"{a}hler-Einstein metric). 
\par For any holomorphic tangent vector $v\in T_o(M)$, we can define the null space of $v$ with respect to the holomorphic bisectional curvature $\mathcal{N}_v=\{w\in T_o(M): R(v,\bar{v}, w, \bar{w})=0\}$, we can show that if $v$ is of maximal rank, $\mathcal{N}_v=0$ (cf.\cite{MR1081948}).  Define $H_v(w, w')=R(v,\bar{v}, w, \bar{w'})$, if $v$ is a characteristic vector of unit length tangent to a minimal rational curve $C$, then $T_o(M)=\mathbb{C}v+\mathcal{H}_v+\mathcal{N}_v$ being decomposed as eigenspaces of $H_v$ where $\mathbb{C}v, \mathcal{H}_v, \mathcal{N}_v$ are corresponding to eigenvalues $2,1,0$ respectively by choosing suitable K\"{a}hler-Einstein metric. Moreover the Grothendieck decomposition of $T(M)$ over $C$ is $T(M)=\mathcal{O}(2)\oplus \mathcal{O}(1)^{\dim\mathcal{H}_v}\oplus \mathcal{O}^{\dim\mathcal{N}_v}$ accordingly.

\par For general knowledge on Hermitian symmetric spaces, we refer \cite{MR1834454}, \cite{MR0404716} and \cite{MR1081948} to the readers.

\subsection{Transitivity of the splitting submanifolds}

The reader may refer the graph given in Section 1.2 in \cite{MR2190340}. We just state this briefly for later discussion. Suppose $A\subset B \subset C$ are complex manifolds, if $A$ is splitting in $B$ and $B$ is splitting in $C$, we have $A$ is splitting in $C$. On the other hand if $A$ is splitting in $C$. then $A$ is splitting in $B$.
This tells us that for example if $A$ is contained in a linear subspace of an irreducible compact Hermitian symmetric space $M$ and $A$ is not a linear subspace, then $A$ is not splitting in $M$, by the result of Van de Ven.

\section{Rational homogeneity of the submanifolds}
\label{Rational-homogeneity}

We prove the following proposition to reduce the problem to a Lie algebra equation.
\begin{proposition}
	Let $M=G/P$ be a rational homogeneous manifold of Picard number one and $S\subset M$ be a submanifold with splitting tangent sequence, $\mathfrak{g}=Lie(G), \mathfrak{p}=Lie(P)$. Fix a gradation for $\mathfrak{g}=\oplus^{m}_{i=-m}\mathfrak{g}_i$ with $T_o(M)\cong \oplus_{i=1}^m\mathfrak{g}_i$ and $\mathfrak{p}=\oplus^{0}_{i=-m}\mathfrak{g}_i$,
	If the following Lie algebra equation has a solution $\eta\in \mathfrak{p}$ for any given  $\zeta\in \oplus^m_{i=1}\mathfrak{g_i}$,
	\[[\eta,\zeta]=\zeta \mod \mathfrak{p}\]
	Then $S$ is rational homogenous.
\end{proposition}
\begin{proof} As we know $S$ can be decomposed as $S\cong A \times N$ by the Theorem of Borel and Remmert (cf.\cite{MR0145557}), where $A$ is an abelian variety and $N$ is a rational homogeneous space. Assume that $A$ is not trivial and we still use $A$ to denote the product $A\times \{pt\}\subset A\times N$, then the tangent sequence $0\longrightarrow T(A) \longrightarrow T(S)|_{A} \longrightarrow N_{A|S}\longrightarrow 0$ splits. By the transitivity property we know the tangent sequence $0\longrightarrow T(A) \longrightarrow T(M)|_A \longrightarrow N_{A|M}\longrightarrow 0$ splits. Without loss of generality we may assume that the reference point $o=eP \in G/P=M$ lies on $A$. We know $\mathfrak{g}$ can be identified with the Lie algebra of global holomorphic vector fields and $\mathfrak{p}$ can be identified with the Lie algebra of global holomorphic vector fields vanishing at the reference point.
	
	Now we consider two global holomorphic vector fields $\eta\in \mathfrak{p}$ and $\zeta\in \mathfrak{g}$ with $\zeta(o)\in T_{o}(A)$. Since $\eta(o)=0$, we have $\nabla_{\eta} \zeta (o)=0 $ for any affine connection $\nabla$. Also, since $\eta(o)=0$, $\nabla_{\zeta} \eta (o)$ is independent of the choice of connections. Therefore we use a torsion-free connection $\nabla$ on $M$, then $(\nabla_{\eta} \zeta - \nabla_{\zeta} \eta)(o)=[\eta,\zeta](o)$. On the other hand, splitting condition gives the restriction and projection map $\pi$ on the space of global holomorphic vector fields, i.e. $\eta |_{A} =\eta_1+\eta_2$, where $\eta_1=\pi(\eta)\in \Gamma(A, T_{A}) $ and $\eta_2 \in \Gamma(A, N_{A|M})$. $\eta(o)=0$ implies that $\eta_{1}(o)=0$, thus $\eta_{1}\equiv 0$ since the tangent bundle of $A$ is trivial. This means that $\eta|_{A}=\eta_{2} \in \Gamma(A, N_{A|M})$, then we have $\nabla_{\zeta}\eta(o) \in N_{{A|M},o} $ since $\eta(o)=0$. This implies $(\nabla_{\eta} \zeta - \nabla_{\zeta} \eta)(o)\in N_{{A|M},o}$. On the other hand, if $\eta \in \mathfrak{p}$ is chosen to be a solution of $[\eta,\zeta]=\zeta \mod \mathfrak{p}$, then $[\eta,\zeta](o)\in T_{o}(A)$. This gives us a plain contradiction. \end{proof}

Since the existence of the solution for the Lie algebra equation is invariant under the isotropy action of $P$, our idea is to classify all the normal forms (representatives) for $P$-orbits in $\mathbb{P}T_o(M)$ and then we show that the equation can be solved for all the representatives by standard linear algebra argument.
For Hermitian symmetric space, the solution can be chosen uniformly to be the central element in $\mathfrak{l}^{\mathbb{C}}$. And the normal forms have been already classified by Proposition \ref{HSS-normal} (it is not necessary in our situation but provides some ideas to deal with more general cases). Moreover we can find a uniform proof for contact case. In fact we can show that the normal form of any vector in the projectivized tangent space can be chosen to be either in $\mathfrak{g}_1$ or $\mathfrak{g}_2$ when $M$ is contact homogeneous.
\begin{theorem}
	Suppose $S\subset M$ is a submanifold with splitting tangent sequence and $M=G/P$ is an irreducible Hermitian symmetric space or a contact homogeneous space of Picard number one, then $S$ must be rational homogeneous. 
\end{theorem}
\begin{proof}
	When $M$ is irreducible Hermitian symmetric, we only have one positive level in the gradation and the central element in $\mathfrak{l}^{\mathbb{C}}$ suffices to give the solution.
	When $M$ is contact, 
	by choosing a Levi factor $L \subset P$ with Lie algebra $\mathfrak{l}^{\mathbb{C}}$, the simple Lie algebra $\mathfrak{g}$ has a graded structure $\mathfrak{g}=\oplus_{k=-2}^{2}\mathfrak{g}_{k}$ with dim$\mathfrak{g}_2=1$, the maximal parabolic subalgebra $\mathfrak{p}=\oplus_{k=-2}^{0}\mathfrak{g}_{k}$ and $\mathfrak{l}^{\mathbb{C}}=\mathfrak{g}_0$. If $\zeta \in \mathfrak{g}_1$, then just choose $\eta$ to be the central element of $\mathfrak{g}_0$ such that $[\eta,v]=kv$ for $v\in \mathfrak{g}_k$, $k=1,2$, especially $[\eta,\zeta]=\zeta$, so we obtain the direct contradiction as in symmetric case. If $\zeta \notin \mathfrak{g}_1$, i.e. $[\zeta]\in \mathbb{P}T_{o}(M)-\mathbb{P}\mathfrak{g}_1$ which is affine, let $\mathcal{U}$ be the unipotent radical of $P$, then the $\mathcal{U}$-orbit (through the adjoint action) of $\mathbb{P}\mathfrak{g}_2$ (identified with a point in the projectivized tangent space) is $\mathbb{P}T_o(M)-\mathbb{P}\mathfrak{g}_1$(cf.Lemma 5 in \cite{MR1681093}). Therefore $\exists \zeta' \in \mathfrak{g}_2$ such that $\zeta'=$ Ad$u.\zeta$ for some $u\in \mathcal{U}$. Thus it suffices to consider two types of normal forms which can be chosen in $\mathfrak{g}_1$ and $\mathfrak{g}_2$  respectively. Then the solution can be found by some $P$-action on the central element in $\mathfrak{g}_0$. 
\end{proof}

The proof for symmetric and contact case is uniform as there are simple properties for the normal forms.
However, we cannot find a uniform proof for all rational homogeneous spaces of Picard number one. Then we have to do case-by-case study to classify the normal forms using certain combinatorial argument, in fact lots of cases can be covered. We have checked the cases whose gradation has depth not greater than $2$.

The classification of normal forms for holomorphic tangent vectors in Hermitian symmetric cases are given according to the maximal strongly orthogonal roots, this can be generalized in some sense for general $G/P$ on the 'abelian part' in the positive part of the gradation  ($\mathfrak{m}^+$ is abelian for Hermitian symmetric spaces), we have 
\begin{proposition}
	For the gradation  $\mathfrak{g}=\oplus_{i=-m}^{m}\mathfrak{g}_i$, if $\mathfrak{g}_k$ is abelian (i.e. $k>\frac{m}{2}$) and if there is a set of maximal strongly orthogonal roots $\Pi=\{\beta_1,\beta_2,...\beta_\ell\}$ with long length in $\mathfrak{g}_k$, then $\{\sum_{j=1}^{\ell} c_je_j, c_j=1 \ \text{or} \ 0\}$ gives normal forms for the $\mathfrak{g}_k$-component under isotropy action of $\exp(\mathfrak{g}_0)$, where $e_j$ is the corresponding normalized root vector of $\beta_j$.
\end{proposition}
\begin{proof}
	Let $\Delta_k$ be the set of all roots in $\mathfrak{g}_k$. We start with $\beta_1$, let $\Pi_1=\{\gamma\in \Delta_k, \beta_1-\gamma \ \text{is a root}\}$, and then let $\Pi_2=\{\gamma\in \Delta_k-\Pi_1,   \beta_2-\gamma \ \text{is a root}\}$, inductively we define $\Pi_j=\{\gamma\in \Delta_k-\sum_{s=1}^{j-1}\Pi_s,   \beta_j-\gamma \ \text{is a root}\}$. Our algorithm can be divided into the following steps
	\begin{description}
		\item[\textbf{Step 1}] For any given vector $\xi$ in $\mathfrak{g}_k$, it can be expressed as a linear combination of root vectors with roots in $\Delta_k$. If the coefficient of $e_1$ is nonzero, we can eliminate all root vectors with roots in $\Pi_1$ using $\mathfrak{g}_{0}$. To see that in this process, when we eliminate one root $\phi\in \Pi_1$ we cannot generate another $\psi\in \Pi_1$ at the same time, we may use the coroot $h_{\beta_1}$. As $\mathfrak{g}_k$ is abelian and $\beta_1$ is a long root, the root spaces are eigenspaces of $ad h_{\beta_1}$ with eigenvalues $2,1,0$ respectively, and more precisely, these three eigenspaces are corresponding to the roots $\{\beta_1\}$, $\Pi_1$ and $\Delta_k-\{\beta_1\}-\Pi_1$. The elimination process decreases the eigenvalues by one or higher, it may generate root vectors with zero eigenvalue but no root vectors with roots in $\Pi_1$ can be generated.
		\item[\textbf{Step 1'}]
		If the coefficient of $e_1$ is zero and there exist some root vectors with nonzero coefficients corresponding to roots in $\Pi_1$, we can generate $e_1$ at first and then start the first step.
		\item[\textbf{Step 2}]
		Now remaining roots are contained in $\Delta_k-\Pi_1$ and the root vector $e_1$ has zero or nonzero coefficient. We then start with $e_2$ and $\Pi_2$, same operations can be done. Also note that the operations on $\Delta_k-\Pi_1-\{\beta_1\}$ will not generate root vectors with roots in $\Pi_1+\{\beta_1\}$, as all root vectors in $\Delta_k-\Pi_1-\{\beta_1\}$ have zero eigenvalues with respect to $ad h_{\beta_1}$ and the operations on $\Delta_k-\Pi_1-\{\beta_1\}$ will not change the eigenvalues so they do not generate root vectors which have been eliminated.
		\item[\textbf{Step 3}] Inductive argument can be used, since $\{\beta_1, \beta_2,..., \beta_\ell\}$ is a set of maximal strongly orthogonal roots, we know under the isotropy action of $\mathfrak{g}_0$ the given $\xi$ can be reduced to $\sum_{j=1}^\ell c_je_j$, and by $\mathbb{C}^*$ action of the Cartan subalgebra, $c_j$ can be reduced to $1$ or $0$, then we finish.  
	\end{description} 
\end{proof} 
\begin{remark}
	We give this proof since there is no Polysphere Theorem and Restricted Root Theorem for general rational homogeneous spaces of Picard number one. These normal forms may be redundant, but it is enough to deal with the problem we consider.
\end{remark}


Next we start to deal with the cases with depth $=2$ in the gradation of $\mathfrak{g}$. Our idea can be simplified as the following. Firstly we obtain a set of maximal strongly orthogonal roots $\Psi_2$ in $\mathfrak{g}_2$ with long length, then we want some \textbf{unique} pairing in $\Delta_1$ in the sense that for any root $\alpha\in\Delta_1$, either it has no pair, or there exists only one root $\beta\in \Delta_1$ such that $\alpha+\beta\in \Psi_2$, this will guarantee that when we use the Ad-action of $\mathfrak{g}_{-1}$ to eliminate the root vectors in $\mathfrak{g}_1$, no new root vectors will be generated except for the root vectors in $\mathfrak{p}=\oplus_{i=-2}^0\mathfrak{g}_i$. Then after enough operations all the remaining roots in $\Delta_1$ must be strongly orthogonal to the remaining roots in $\Delta_2$. Then we can use the following lemma.
\begin{lemma}
	\rm Suppose we have $R_1, R_2$ to be the set of remaining roots in $\Delta_1$ and $\Delta_2$ respectively after some adjoint actions by $\mathfrak{g}_0$ and $\mathfrak{g}_{-1}$. From the above argument we know for $\alpha\in R_1$ and $\beta \in R_2$, they are strongly orthogonal. And all roots in $R_2$ are mutually strongly orthogonal. Then we can solve the equation $[\eta,\zeta]=\zeta\mod \mathfrak{p}$ for any given $\zeta\in \mathfrak{g}_1\oplus \mathfrak{g}_2$, where $\eta\in \mathfrak{p}$.
\end{lemma}
\begin{proof} Let $h=\sum_{\alpha\in R_2}h_{\alpha}$ be the sum of coroots and after $P$-action we write $\zeta=v_1+v_2$, where $v_1$ and $v_2$ are linear combinations of root vectors corresponding to roots in $R_1$ and $R_2$ respectively. Then $[h,v_1]=0$ and $[h, v_2]=2v_2$. Then root vectors with roots in $R_1$ and $R_2$ have eigenvalues $(0,2)$ for $adh$. On the other hand we know $\mathfrak{g}_1$ and $\mathfrak{g}_2$ have eigenvalues $(1,2)$ for $ad\eta_0$ where $\eta_0$ is the central element of $\mathfrak{g}_0$. Then we can solve the equation by combining $h$ and $\eta_0$ linearly. 
\end{proof}

Then the rest work is the combinatorial checking for root systems, i.e. to write down a set of maximal strongly orthogonal roots with long length in $\Delta_2$ and check the unique pairing property. As an example we will give the discussion for $(D_5, \alpha_3)$, we adapt the root system setting as in \cite{MR1274961}.

\begin{example}
	In the case $(D_5, \alpha_3)$, we have $\text{dim}(\mathfrak{g}_2)=3, \text{dim}(\mathfrak{g}_1)=12, \text{dim}(\mathfrak{g}_0)=15$. For completeness, we list the corresponding roots as follows:
	\par $\Delta_0$: $\pm\alpha_1, \pm\alpha_2, \pm\alpha_4, \pm\alpha_5,\pm(\alpha_1+\alpha_2)$
	\par $\Delta_1$: $\alpha_3,
	\alpha_1+\alpha_2+\alpha_3,  \alpha_2+\alpha_3, \alpha_3+\alpha_4, \alpha_3+\alpha_5, \alpha_1+\alpha_2+\alpha_3+\alpha_4, \alpha_2+\alpha_3+\alpha_4, \alpha_1+\alpha_2+\alpha_3+\alpha_5, \alpha_2+\alpha_3+\alpha_5, \alpha_3+\alpha_4+\alpha_5,
	\alpha_2+\alpha_3+\alpha_4+\alpha_5, \alpha_1+\alpha_2+\alpha_3+\alpha_4+\alpha_5$
	\par $\Delta_2$: $\alpha_1+2\alpha_2+2\alpha_3+\alpha_4+\alpha_5, \alpha_1+\alpha_2+2\alpha_3+\alpha_4+\alpha_5, \alpha_2+2\alpha_3+\alpha_4+\alpha_5 $

	\begin{description}
		\item[\textbf{Step 1}] We start from $\mathfrak{g}_2$-component, use Ad-action by $\mathfrak{g}_0$, we remain to have the highest root and  $\Psi_2=\{\alpha_1+2\alpha_2+2\alpha_3+\alpha_4+\alpha_5\} \subset \Delta_2$.
		\item[\textbf{Step 2}]  Since there is at most one root left in $\Delta_2$, the uniqueness for pairing in $\mathfrak{g}_1$ is automatic.
	\end{description}
\end{example}

The Lie algebra equation is actually a system of linear equations if we express it in terms of root vectors. There is a uniform algorithm for cases with depth $\leq 2$ including Hermitian symmetric case as discussed above. For cases with depth $\geq 3$, our strategy is to eliminate root vectors in $\zeta$ as many as possible so that the remaining equation is easy enough to solve by basic linear algebra, but the checking becomes too lengthy to finish, some new ideas may be needed.

The full checking for cases in Theorem \ref{rational splitting} is given in the appendix. Type $A$ are all Grassmannians and for Type $G_2$, $(G_2,\alpha_2)$ is contact homogeneous and $(G_2,\alpha_1)$ is isomorphic to a hyperquadric. We remain to check Type $B,C,D,E,F$ when the gradation has depth $=2$. 

\section{Hermitian symmetry of the submanifolds}
\label{Hermitian-symmetry}

In this part we will furtherly prove that when the ambient space is irreducible Hermitian symmetric, the compact splitting submanifolds are Hermitian symmetric as well. We still make use of the restriction and projection map $\pi$ on the global vector fields of the ambient space.
\par More precisely we will use the special property of Euler vector field (which is vanishing to the first order at the reference point) that the restriction of Euler vector field to any affine subspace is a tangent Euler vector field of the subspace. The key point is that the central element in $\mathfrak{g}_0$ gives an Euler vector field if and only if the space is symmetric. Contradiction will arise if the splitting submanifold is not symmetric since the map $\pi$ on the Euler vector field on the ambient space will give an Euler-like vector field on the subspace (up to a higher order error term).

We firstly give a well-known fact about the Euler vector field,

\begin{lemma}
	If $(z_1,z_2,...z_n)=F(w_1,...,w_n)=(f^1(w_1,...,w_n),...,f^n(w_1,...,w_n))$ is a change of complex coordinate with $F(0,...,0)=(0,...,0)$ (i.e. the map $F$ is a biholomorphism preserving the origin), then the Euler vector field \[\sum_{i=1}^nw_i\frac{\partial}{\partial w_i}=\sum_{i=1}^nz_i\frac{\partial}{\partial z_i}+\text{term of higher vanishing order}\geq 2\], i.e. under change of complex coordinate the first order part in the vector field is also the Euler vector field in the new coordinate.
\end{lemma}
\begin{proof}
	The proof is given by elementary checking. Write the inverse map as $(w_1,w_2,...w_n)=F^{-1}(z_1,...,z_n)=(g^1(z_1,...,z_n),...,g^n(z_1,...,z_n))$. We can compute 
	\[
	\begin{aligned}
	\sum_{i=1}^nz_i\frac{\partial}{\partial z_i}&=\sum_{i,j=1}^nf^i\frac{\partial g^j}{\partial z_i}\frac{\partial}{\partial w_j}\\
	&=\sum_{i,j,k=1}^nw_k\frac{\partial f^i}{\partial w_k}\frac{\partial g^j}{\partial z_i}\frac{\partial}{\partial w_j}+\text{term of higher vanishing order}\geq 2\\
	&=\sum_{j,k=1}^nw_k\delta^j_k\frac{\partial}{\partial w_j}+\text{term of higher vanishing order}\geq 2\\
	&=\sum_{j=1}^nw_j\frac{\partial}{\partial w_j}+\text{term of higher vanishing order}\geq 2
	\end{aligned}
	\] 
\end{proof}

\par Then we start to prove that $S$ is Hermitian symmetric.
\begin{proof}[Proof of the Hermitian symmetry of $S$]
	Let $o\in S\subset M$ be the reference point, $n=\dim(M), m=\dim(S)$.
	Let $c\in \mathfrak{l}^\mathbb{C}$ be the central element which gives a global vector field on $M$ vanishing to the first order at $o$, and we can choose coordinate (may be different from the Harish-Chandra coordinate) such that the germ of the submanifold $S$ is affine in this coordinate. We know $c$ is the Euler vector field under the Harish-Chandra coordinate, then from the lemma, written in the chosen coordinate $c$ can be represented as
	\[c=\sum_{i=1}^nz_i\frac{\partial}{\partial z_i}+\text{term of higher vanishing order}\geq 2\]
	whose first order part is the Euler vector field. Without loss of generality the germ can be written in this coordinate as $z_{m+1}=\cdots =z_n=0$ and then after the restriction and projection map we know the global vector field $\pi(c)$ on $S$ written in this local coordinate is 
	\[\pi(c)=\sum_{i=1}^mz_i\frac{\partial}{\partial z_i}+\text{term of higher vanishing order}\geq 2\]
	On the other hand we consider the global  vector fields on $S$ vanishing at $o$. Since $S$ is a rational homogeneous space, we write it as $S=\Pi_{k}G_k'/P_k'$ for connected simple Lie groups $G'_k$ and corresponding parabolic subgroups $P_k'$. Since each factor in $S$ is also splitting in $M$, without loss of generality we consider the case when $S$ has only one factor. 
	
	\par In this case we write $S=G'/P'$. The Lie algebra of $P'$ is denoted by $\mathfrak{p}'$ which is also the Lie algebra of the global vector fields on $S$ vanishing at $o$. For any $\eta\in \mathfrak{p}'$, it can be expressed in the coordinate as 
	\[ \eta=\sum_{j=1}^m \eta_j(z_1,...,z_m)\frac{\partial }{\partial z_j}   \] where $\eta_j(z_1,...,z_m)$ are holomorphic and vanish at $o$. We now compute $[\pi(c),\eta]$ where 
	\[\pi(c)\eta=\sum_{i,j=1}^mz_i\frac{\partial\eta_j}{\partial z_i}\frac{\partial }{\partial z_j}+\text{term of higher vanishing order}\geq 2\]
	and 
	\[\eta\pi(c)=\sum_{j=1}^m\eta_j\frac{\partial}{\partial z_j}+\text{term of higher vanishing order}\geq 2\]
	Taking the Taylor expansion on $\pi(c)\eta$ and $\eta\pi(c)$ at $o$ we know the first order terms are the same and hence $[\pi(c),\eta]$ is a global vector field on $S$ vanishing to the second order at $o$. 
	\par Fix a graded structure on $\mathfrak{g}'=Lie(G')$ with $\mathfrak{p}'=\oplus_{k=-q}^{0}\mathfrak{g}'_k$ where $q$ is the depth of the gradation, we know $\mathfrak{g}'_{-q}$ is identified with the Lie algebra of global vector fields vanishing to the second order. The argument above implies that there exists a nonzero element $\pi(c)\in \mathfrak{p}'$ such that $[\pi(c),\mathfrak{p}']\subset \mathfrak{g}'_{-q}$, and hence there exists a nonzero element in $\mathfrak{g}_0'$ which centralizes $\oplus_{k=-q+1}^{0}\mathfrak{g}'_k$, this is possible if and only if $q=1$ (cf.\cite[Lemma 3.2]{MR1274961}), i.e. $S$ is Hermitian symmetric.  
\end{proof}

\section{A differential geometric proof for the classification of splitting submanifolds in the hyperquadric: an example for a possible new approach}
\label{new-method}

We use the map $\pi$ on the Euler vector field in the last section and it is a global holomorphic vector field on $M$ vanishing to the first order at $o$. A natural idea is to use the map $\pi$ on the global holomorphic vector fields on $M$ vanishing to the second order at $o$ to investigate more information on the splitting submanifolds.
From this we may obtain a possible procedure to find the splitting submanifolds for explicit examples by classifying the tangent subspaces and using some integrability argument. This may be a possible new approach for the classification, at least when the ambient space is an irreducible compact Hermitian symmetric space of rank two. We give a new proof for Jahnke' result (cf. Theorem \ref{thm1}) in this part as an explanation.
\subsection{Degree of the embedding}
\label{Degree-embedding}
From previous section we know $S\subset M$ is Hermitian symmetric when $M$ is irreducible Hermitian symmetric but we know nothing about the embedding. In this part we show the following proposition on the degree of the embedding of $S$ in $M$ when $M$ is of rank two. Here minimal rational curves for product of irreducible compact Hermitian symmetric spaces means the pullback of minimal rational curves for each factor (i.e. rational curves with muti-degree (0,...,1...,0)).

\begin{proposition}\label{linearity rank-2}
	Let $M$ be an irreducible rank two Hermitian symmetric space of compact type. If $f:S\hookrightarrow M$ is a splitting submanifold other than projective spaces, then the embedding $f$ is linear, i.e. $f$ preserves the minimal rational curves.
\end{proposition}

We firstly give a theorem proven in \cite{MR3647134} which will be used in our proof.
\begin{theorem}\label{dual char}\rm(\textbf{Mok-Ng, \cite[p.10]{MR3647134}})
	Suppose $M$ is an irreducible Hermitian symmetric space equipped with some choice of canonical K\"{a}hler-Einstein metric,  and $S$ is a compact characteristic complex submanifold with splitting tangent sequence in $M$. Then $(S,g|S)\hookrightarrow(M,g)$ is totally geodesic and moreover $S$ is invariantly totally geodesic and it must be a projective space. Here characteristic complex submanifold means that all tangent vectors of the submanifold are characteristic.
\end{theorem}

Then we start to prove the proposition.
\begin{proof} We know $S$ is Hermitian symmetric and it is not a projective space. If $S$ is reducible, i.e. the Picard number is at least two, then minimal rational curves of $S$ are defined by pullback of the minimal rational curves of each factor in $S$. We choose a minimal rational curve $C$ in $S$, then $C$ is splitting in $S$. Since $S$ is splitting in $M$ and the Grothendieck decomposition of $T(S)$ over $C$ has at least one $\mathcal{O}$-factor we know there is also an $\mathcal{O}$-factor in the Grothendieck decomposition of $T(M)$ over $C$. We claim that tangent vectors at every points on $C$ are characteristic. For any point $x\in C$, we denote the tangent vector at $x$ by $e_{\alpha}$. We choose some canonical K\"{a}hler-Einstein metric $g$ on $M$ with nonnegative bisectional curvature. Then there is a Hermitian metric $g^*$ on the cotangent bundle $T^*_{M}$ such that $T^*_{M}$ is of non-positive curvature. Let $\eta$ be a nowhere zero holomorphic section of $L=\mathcal{O}$ (as a subbundle of $T^*_{M}$) with unit length. Then
	\[0=Ric(L,g^*|L)(\alpha, \bar{\alpha})\leq R^{L}_{\alpha\bar{\alpha}\eta\bar{\eta}}\leq R^{T^*_{M}}_{\alpha\bar{\alpha}\eta\bar{\eta}}\leq0\] which implies $0=R^{L}_{\alpha\bar{\alpha}\eta\bar{\eta}}= R^{T^*_{M}}_{\alpha\bar{\alpha}\eta\bar{\eta}}$. Then for the tangent vector $e_\alpha$, nullity for the bisectional curvature on $M$ is nonzero. Also since $M$ is of rank two, each tangent vector is either characteristic or of rank two. But we know for a rank two vector, the nullity is zero. Thus $e_{\alpha}$ is a characteristic vector and all tangent vectors of $C$ are characteristic on $M$. By Theorem \ref{dual char}, we know $C$ is totally geodesic (actually invariantly totally geodeisc). Hence $C$ is a minimal rational curve in $M$, i.e. $f$ is a linear embedding.  \end{proof}
\begin{remark}
	When $M$ is of rank two, from the proof above we know if a rational curve $C\subset M$ is splitting and the Grothendieck decomposition of $T_M$ over $C$ has $\mathcal{O}$-factors, then every tangent vector of $C$ is characteristic thus $C$ is invariantly totally geodesic, i.e. $C$ is a minimal rational curve in $M$.
\end{remark}

\begin{corollary}\label{projective space in rank-2}
	Let $M$ be an irreducible rank two compact Hermitian symmetric space. If
	$S=X_1\times X_2 \times \cdot\cdot\cdot \times X_m (m\geq 2)$ is a product of irreducible compact Hermitian symmetric spaces which is splitting in $M$, then there is no factor $X_i$ with rank$(X_i)\geq 2$ and furthermore $S$ must be a product of linear projective spaces with only two factors.
\end{corollary}
\begin{proof} If a product of irreducible compact Hermitian symmetric space $S=X_1\times X_2 \times \cdot\cdot\cdot \times X_m$ is splitting, then from the previous proposition each factor $X_i$ is linear, i.e. $\mathcal{O}_{M}(1)|_S=\sum_{1\leq i\leq m}\mathcal{O}_{X_i}(1)$. And by polysphere theorem the question can be reduced to prove that there is no $\mathbb{P}^1\times Q^2=\mathbb{P}^1\times \mathbb{P}^1\times \mathbb{P}^1$ splitting in $M$. If $Y=\mathbb{P}^1\times Q^2$ is splitting in $M$, choose any smooth rational curve $C$ in $Q^2$, we know it is splitting in $Q^2$ (cf. \cite[Proposition 4.3]{MR2190340}) thus the pullback of $C$ in $Y$ (which is also denoted by $C$) is splitting in $Y$ and hence in $M$. Then we have $N_{C|M}=N_{C|Y}\oplus N_{Y|M}|_{C}$. Since $N_{C|Y}$ has at least one $\mathcal{O}$-factor, we can follow the remark above to obtain that $C$ is a minimal rational curve in $M$. But we can choose $C$ in $Q^2$ of any degree, a contradiction. Hence there is no factor $X_i$ with rank $\geq 2$. This also rule out the case that $S$ is a product of at least three linear projective spaces since we can always find a $\mathbb{P}^1\times Q^2$ splitting in such $S$. So $S$ must be a product of two linear projective spaces. 
\end{proof}

\subsection{The restriction on the tangent subspaces}
\label{Restriction-tangent}
In this part we give a restriction on the tangent subspaces by the map $\pi$ on global holomorphic vector fields 
vanishing to the second order at $o$ on the ambient space. The restriction on the tangent subspaces will provide restriction on the embedding of $S$ in further discussion.

We know the following fact 
\begin{proposition}[cf.\cite{MR2178704}]\label{characterization by vector fields}
	Let $d_2(x)$ denote the dimension of space of global holomorphic vector fields vanishing to the second order at $x$ on a rational homogeneous manifold $S$, then $d_2(x)\leq \dim(S)$ and the bound is sharp if and only if $S$ is Hermitian symmetric.
\end{proposition}

Let $\mathfrak{f}_M, \mathfrak{f}_S$ denote the space of global holomorphic vector fields vanishing to the second order at $o$ on $M$ and $S$ respectively.
For any $Z\in \mathfrak{f}_{M}$, its 2-jet at $o$ gives an element in $S^2T^*_{o}(M)\otimes T_{o}(M)$. We know there is an embedding \[j_{M}:\mathfrak{f}_{M}\cong\mathfrak{m}^{-}\hookrightarrow S^2T^*_{o}(M)\otimes T_{o}(M)\cong Hom(S^2\mathfrak{m}^{+}, \mathfrak{m}^+)\] where $Z$ is corresponding to some $\bar{w}\in \mathfrak{m}^{-}$ (here and henceforth the complex conjugate is with respect to the compact real form of $\mathfrak{g}$) and $j_{M}(\bar{w})$ can be written as $[-,[\bar{w},-]]$. If we choose some basis $\{v_1,v_2,...,v_{n}\}$ for $T_{o}(M)$, then $[-,[\bar{w},-]]=\sum_{i,j=1}^n[v_i,[\bar{w},v_j]v^*_i\otimes v^*_j$. Similarly we have the embedding 
\[j_{S}:\mathfrak{f}_{S} \hookrightarrow S^2T^*_{o}(S)\otimes T_{o}(S) \]

For simplicity we use the pairing $(V,W)$ to denote the candidates for the tangent subspace $T_o(S)$ and its complement.
The decomposition $T_o(M)=V\oplus W$ induces a decomposition of cotangent space $T^*_o(M)=Ann(W)\oplus Ann(V)$, where $Ann(W)\cong T^*_{o}(S)$ (from the bundle isomorphism). We have a projection $\pi_o$ from $S^2T^*_{o}(M)\otimes T_{o}(M)$ to $S^2T^*_{o}(S)\otimes T_{o}(S)$ induced from the decompositions and hence 
have a projection of holomorphic 2-jets.
We can transform the argument on projection of vector fields to the argument on projection of 2-jets from the following lemma. Then the problem can be 'linearized' to the tangent space at one point. Without ambiguity $\pi_o$ also denotes the projection on the tangent and cotangent spaces when we consider the projection of 2-jets.

\begin{lemma}
	The following diagram is commutative on $\ker(\pi)\subset \mathfrak{f}_M$
	\[	\xymatrix{
		\mathfrak{f}_{M} \ar[rr]^-{j_{M}}\ar[d]_{\pi} & & S^2T^*_{o}(M)\otimes T_{o}(M) \ar[d]^{\pi_o} \\
		\mathfrak{f}_{S}\ar[rr]^-{j_{
				S}}	&  & S^2T^*_{o}(S)\otimes T_{o}(S)
	}\]
	and hence $\dim \ker(\pi)\leq \dim \ker(\pi_o|_{Im(j_M)})$.
	
\end{lemma}

\begin{proof}
	Choose some element $\overline{\xi}\in \mathfrak{m}^-$ and write in Harish-Chandra coordinate $\overline{\xi}=\sum_{j,k,\ell}A_{k\ell}^jz_kz_{\ell}\frac{\partial}{\partial z_{j}}$, with $A^j_{k\ell}=A^j_{\ell k}$. Then $j_{M}(\overline{\xi})=\sum_{j,k,\ell}A_{k\ell}^jdz_k\odot dz_{\ell}\otimes \frac{\partial}{\partial z_{j}}$ where $\odot$ denotes the symmetric tensor. Locally a germ of the submanifold $S$ can be written as the image of some holomorphic map $\tau:(s_1,...,s_m)\rightarrow (f^1(s_1,...,s_m),...,f^n(s_1,...,s_m))$ from $\mathbb{C}^m$ to $M$ with $\tau(0,...,0)=(0,...,0)=o$, then the restriction of $\overline{\xi}$ on $S$ is $\sum_{j,k,\ell}A_{k\ell}^jf^kf^{\ell}\frac{\partial}{\partial z_{j}}$. 
	\par It suffices to check that $\pi(\overline{\xi})=0$ implies $\pi_o\circ j_{M}(\overline{\xi})=0$. $\pi(\overline{\xi})=0$ gives that $\sum_{j,k,\ell}A_{k\ell}^jf^kf^{\ell}\frac{\partial}{\partial z_{j}}\in N_{S|M,x}$ for $x$ in an open neighborhood of $o$. Also since $\overline{\xi}$ vanishes at $o$ to the second order we know \[\nabla_{\frac{\partial}{\partial s_a}}\nabla_{\frac{\partial}{\partial s_b}}\sum_{j,k,\ell}A_{k\ell}^jf^kf^{\ell}\frac{\partial}{\partial z_{j}}|_{(0,...,0)}\in W=N_{S|M,o}\]
	This gives that \[\sum_{j,k,\ell}A_{k\ell}^j\frac{\partial f^k}{\partial s_a}\frac{\partial f^{\ell}}{\partial s_b}\frac{\partial}{\partial z_{j}}|_{(0,...,0)}\in W\]
	for any $1\leq a,b \leq m$. Hence \[\pi_o\circ j_{M}(\overline{\xi})=\sum_{a,b}\sum_{j,k,\ell}A_{k\ell}^j\frac{\partial f^k}{\partial s_a} \frac{\partial f^{\ell}}{\partial s_b}|_{(0,...,0)}ds_a\odot ds_b\otimes \pi_o(\frac{\partial}{\partial z_{j}})=0\]
	here we also use $\pi_o$ to denote the projection from $T_o(M)$ to $T_o(S)$.
	
\end{proof}

From the rational homogeneity of $S$ we know $\dim\ker(\pi)\geq n-m$. Also from the lemma we know if a subspace $V\subset T_o(M)$ satisfies
$\dim \ker(\pi_o|_{Im(j_M)})<n-m$ for any choice of complementary subspace $W \subset T_o(M)$, then $V$ can not be a candidate for a tangent space of $S$,
here  $\ker(\pi_o|_{Im(j_M)})=\{\overline{\xi}\in \mathfrak{m}^-: [v,[\overline{\xi},v]]\in W, \forall v\in V \}$. This provides some information on the tangent subspaces and then the embedding of $S$, at least we can compute some examples. Following the idea we will give a differential geometric proof for the classification of splitting submanifolds with $\dim\geq 2$ in the hyperquadric $Q^n$. The discussion for the case $M=G(2,3)$ is also considered in the author's thesis \cite{DingThesis}, which will not appear in this article as the computation is more tedious.
\begin{remark}
	We originally want to follow the idea to show that $\dim(\ker (\pi))\leq \dim\ker(\pi_o|_{Im(j_M}))\leq n-m$ for any choice of $V,W$ and together with the fact that $\dim\ker(\pi)\geq n-m$ to show the Hermitian symmetry of $S$. However the inequality $\dim\ker(\pi_o|_{Im(j_M}))\leq n-m$ is not easy to check in general although we can find that it holds for some examples. So we use the Euler vector field instead in previous section and keep this idea in the discussion on concrete examples.
\end{remark}

At the end of this part we give a simple lemma on the $P$-invariance of $\dim\ker(\pi_o|_{Im(j_M)})$.
\begin{lemma}\label{invariant P-action}
	The dimension of $\ker(\pi_o|_{Im(j_M)})$ is invariant under $P$-action, i.e. for any $p\in P$, $\dim\ker(\pi_o|_{Im(j_M)})$ is the same for the decomposition $T_{o}(M)=V\oplus W$ and $T_{o}(M)=Adp.V\oplus Adp.W$.
\end{lemma}
\begin{proof} If the decomposition at $o$ is given by $T_{o}(M)=V\oplus W$, we know \[\begin{aligned}\ker(\pi_{o}|_{Im(j_M}))&=\{\overline{\xi}\in \mathfrak{m}^-: [v,[\overline{\xi},v]]\in W, \forall v\in V \}\\&=\{\overline{\xi}\in\mathfrak{m}^-: R(v_1,\overline{\xi}, v_2, \overline{u})=0, \forall v_1,v_2\in V, u\in W^\bot)\}
	\end{aligned}
	\]
	Hence $\ker(\pi_{o}|_{Im(j_M)})=\overline{([[V, \overline{W^\bot}],V])^\bot}$ where the orthogonal complement is with respect to some invariant metric induced from the Killing form. Since $\dim([[V, \overline{W^\bot}],V])$ is invariant under $P$-action, so is $\dim\ker(\pi_{o}|_{Im(j_M)})$.
\end{proof}

\subsection{Proof for the classification by projection of 2-jets}
\label{hyperquadric-two-jets}

Motivated by the method of restriction and projection on global vector fields, we will give a new proof for the classification for compact splitting submanifolds $S$ with dimension $m \geq 2$ in the hyperquadric $Q^n$ based on differential geometry, which is different from the algebraic geometric argument given by Jahnke\cite{MR2190340} based on the fact that the quadrics are only projective manifolds $X$ such that $\bigwedge^2 T_X$ is ample. 

We use $(V,W)$ to denote the candidates for the tangent subspace and its complement, we will find that $\dim(\ker(\pi_o|_{Im(j_M)}))\leq n-m$ for any pair $(V,W)$ and $\dim(\ker(\pi_o|_{Im(j_M)}))=n-m$ only if the tangent subspaces are of some 'standard types'. Then from the classification on the tangent spaces of $S$, together with some integrability argument we obtain the classification for $S$. 

For later discussion, we give a description for the holomorphic conformal structure on $Q^n$ (which is actually the projective second fundamental form of $Q^n$). This is an Aut$(Q^n)$-invariant holomorphic section of $\Gamma(Q^n, S^2T^*(Q^n)\times \mathcal{O}(2))$ and 
on the reference point $o$, identifying the fiber of $\mathcal{O}(2)$ with $\mathbb{C}$ by some linear isomorphism, it can be written as $q_o=e_1^*\odot e_1^*+\cdots+e_n^*\odot e_n^*$ where $\{e_1,...,e_n\}$ is chosen as an orthonormal basis of $T_o(Q^n)$ consisting of rank two vectors. One can see the VMRT of $Q^n$ at $o$ is actually defined by the kernel of $q_o$.

We then consider the compact splitting submanifolds $S\subset Q^n$ with $\dim(S)=m\geq 2$.

\subsubsection{Computation of 2-jets for $Q^n$} Choose $o$ on $S$ as a reference point of $Q^n$ and in Harish-Chandra coordinate any nonzero holomorphic vector field on $Q^n$ vanishing to the second order at $o$ induces a 2-jet as follows, here $\odot$ denotes the symmetric tensor on $S^{2}T^*_o(Q^n)$ and $\frac{\partial}{\partial z_i}=e_i$. \[\sum_{1\leq i,j,k\leq n} A^{k}_{ij}dz^i\odot dz^j\otimes \frac{\partial}{\partial z_k}\]
which is in $Hom(S^2T_o(Q^n), T_o(Q^n))\cap \mathfrak{l}^{\mathbb{C}}\otimes T^*_o(Q^n)$ and $A$ is chosen to satisfy $A^k_{ij}=A^k_{ji}$. In this case $(A_i)^k_j\in \mathfrak{so}(n,\mathbb{C})\oplus \mathbb{C}id$, so we know if $i,j,k$ are different from each other, $A^{k}_{ij}=-A^{j}_{ik}=A^{i}_{kj}=-A^{k}_{ji}$ and then $A^{k}_{ij}=0$. All other coefficients are determined by $A_{ik}^{k}$. So any such 2-jet can be written as 
\[Z=\sum_{j=1}^n a_j(D\otimes \frac{\partial}{\partial z_j}+dz^j\odot B)\]
where $D=-\sum_{s=1}^n (dz^s)^2, B=2\sum _{s=1}^ndz^s\otimes \frac{\partial}{\partial z_s}$ and $a_j\in \mathbb{C}$.

\subsubsection{Normal forms (representatives under $P$($K^{\mathbb{C}}$)-action) for tangent subspaces}
For a subspace $V\subset T_o(Q^n)$ with $\dim V=m$, according to the rank of holomorphic conformal structure $q_o|_{V}$, the normal forms of $m$-planes in $T_o(Q^n)$ under $K^{\mathbb{C}}$-action (representatives of $K^{\mathbb{C}}$-orbit on $Gr(m,T_o(Q^n))$) can be classified as \[\mathbb{C}\{\frac{\partial}{\partial z_1},...,\frac{\partial}{\partial z_k}, \frac{\partial}{\partial z_{k+1}}+i\frac{\partial}{\partial z_{m+1}},..., \frac{\partial}{\partial z_m}+i\frac{\partial}{\partial z_{2m-k}}\}\]
where $k$ is the rank of $q_{o}|_V$.
Assume that $k\neq 0, k<m$(i.e. $V$ is neither  characteristic nor nondegenerate with respect to $q_o$).

\subsubsection{Projection of 2-jets}
By Lemma \ref{invariant P-action} we know it suffices to compute for the normal forms of $V$.
The splitting condition gives the projection on the space of 2-jets 
\[\pi_o(Z)=\pi_o(D)\otimes \pi_o(\sum_{j=1}^na_j\frac{\partial}{\partial z_{j}})
+\pi_o(\sum_{j=1}^ka_jdz^j+\sum_{j=k+1}^m(a_j+ia_{j+m-k})dz^j)\odot  \pi_o(B)	\]
since we know \[\ker(\pi_{o,T_o^*(Q^n)})=Ann(V)=\mathbb{C}\{dz^{2m-k+1},...,dz^n,dz^{k+1}+idz^{m+1},...,dz^m+idz^{2m-k}\}\]
\par Also we have 
\[\pi_o(B)=2\sum_{s=1}^k\pi_o(dz^s)\otimes \frac{\partial }{\partial z_s}+2\sum_{s=k+1}^{m}\pi_o(dz^s)\otimes (\frac{\partial }{\partial z_s}+i\frac{\partial }{\partial z_{s+m-k}})\]

Here we note that $\pi_o(B)$ is not decomposable and $\pi_o(B)\neq 0, \pi_o(D)\neq 0$ (since $k\neq 0$ and $k<m$), so if $\pi_o(Z)=0$ we have \[\begin{cases}
\pi_o(\sum_{j=1}^na_j\frac{\partial}{\partial z_{j}})=0 \\
\sum_{j=1}^ka_j\pi_o(dz^j)+\sum_{j=k+1}^m(a_j+ia_{j+m-k})\pi_o(dz^j)=0
\end{cases}
\]
The first condition is precisely coincident with $\ker(\pi_{o,T_o(Q^n)})$. For the second condition, since $\pi_o(dz^j)(1\leq j\leq m)$ are linearly independent we know $a_j=0$ for $1\leq j\leq k$ and $a_j+ia_{j+m-k}=0$ for $k+1\leq j\leq m$. Since $k<m$ and $ \pi_o(\frac{\partial}{\partial z_{j}}+i\frac{\partial}{\partial z_{j+m-k}})\neq 0$ for $k+1\leq j\leq m$, we know the second condition gives an $(n-m)$-dimensional vector space which is not 
coincident with $\ker(\pi_{o,T_o(Q^n)})$.
This implies that the dimension of $
\ker(\pi_o|_{Im(j_{M})})$ is strictly less that $n-m$. Hence we have
\begin{proposition}\label{quadric normal form}
	If $S$ is a compact splitting submanifold in $Q^n$ with dimension at least 2, then for any $x\in S$, the tangent space $T_x(S)$ is either characteristic or nondegenerate with respect to the holomorphic conformal structure of $Q^n$.
\end{proposition}

Furthermore we can prove

\begin{corollary}\label{complementary}
	
	If $T_o(S)$ is nondegenerate with respect to the holomorphic conformal structure of $Q^n$, then the complementary subspace $N_{o,S|Q^n}$ in the direct sum is uniquely determined and it is the annihilator
	of $T_o(S)$ with respect to the holomorphic conformal structure.
\end{corollary}
\begin{proof}
	Let $w_i=\pi_o(\frac{\partial}{\partial z_i})$ and $w^i=\pi_o(dz^i)$.
	Inheriting the computation above we know if $k=m\geq 2$, 
	\[
	\pi_o(Z)=2\sum_{j=1}^ma_jw^j \odot (\sum_{s=1}^m w^s\otimes w_s)-(\sum_{j=1}^mw^j\odot w^j)\otimes (\sum_{s=1}^na_sw_s)	\]
	Then $\pi_o(Z)=0$ implies \[
	\begin{cases}
	a_1=\cdots=a_m=0\\
	\sum_{s=m+1}^na_sw_s=0
	\end{cases} \]
	Thus the kernel of $\pi_o$ is $(n-m)$-dimensional if and only if $w_{m+1}=\cdots=w_n=0$, i.e. $N_{o,S|Q^n}$ is uniquely determined and it is the annihilator
	of $T_o(S)$ with respect to the holomorphic conformal structure.   
\end{proof}

Also we give a theorem proven in \cite{Zhang2015} related to the integrability.
\begin{theorem}\label{second order quadric}\rm(\textbf{Zhang, \cite[Main Theorem]{Zhang2015}}) If for any $x\in S\subset Q^n$, there is a smooth linear section subquadric tangent to $S$ at $x$, and germ of any minimal rational curve $L$ of $Q^n$ passing through $x$
	with $T_x(L)\subset T_x(S)$ lies on $S$, then $S$ itself is a smooth linear section subquadric.
\end{theorem}

Then we can give a new proof on the theorem obtained by Jahnke\cite{MR2190340}.

\begin{theorem}
	Any splitting submanifold with dimension $m\geq2$ in the hyperquadric $Q^{n}$ is a linear subspace or a smooth complete linear section subquadric. Any splitting curve in $Q_{n}$ is rational.
\end{theorem}
\begin{proof} 		For an $m$-dimensional splitting submanifold $S\subset M=Q^n$ with $m\geq 2$, the set of points with characteristic tangent spaces is closed. We claim that it is also open so that tangent spaces of $S$ are all the same isomorphism type. If $\{x_k\}_k$ is a sequence of points on $S$ with tangent spaces tangent to some smooth linear section subquadrics, from Corollary \ref{complementary} we know complementary subspace of $T_{x_k}(S)$ must be its annihilator with respect to the holomorphic conformal structure. Suppose $x_k$ goes to $x$. Then the complementary subspace of $T_x(S)$ must be its annihilator with respect to the holomorphic conformal structure as well, which implies that $T_x(S)$ is also tangent to a smooth linear section subquadric. Thus the claim holds.
	
	Then we firstly consider the case when every tangent space of $S$ is nondegenerate with respect to the holomorphic conformal structure of $Q^n$ (i.e. $S$ is every point tangent to some smooth linear section subquadrics). If $S$ is not a projective space, then we know the embedding is linear from Proposition \ref{linearity rank-2} and hence it satisfies the condition in Theorem \ref{second order quadric}. Therefore it must be a smooth linear section subquadric.
	If $S$ is a projective space in $Q^n$ whose tangent space is nondegenerate with respect to the holomorphic conformal structure of $Q^n$, then $S$ is not a degree one subspace and hence there exists a high-degree projective space $\mathbb{P}^2$ which is splitting in $Q^n$. And this projective space is tangent to some smooth linear section subquadrics isomorphic to $Q^2$ at every point. Then there exists a holomorphic section for $\mathbb{P}T(\mathbb{P}^2)$ which gives a proper subbundle for $T(\mathbb{P}^2)$, a contradiction.
	
	If $S$ is every point characteristic, then by Theorem \ref{dual char}, $S$ is a linear subspace in $Q^n$. 
\end{proof}

\appendix
\section{Root system checking for cases with depth $\leq 2$}
\label{depth2}
We adapt the root system notations from \cite{MR1274961} for $B,C,D$ types and \cite{MR0240238} for $E,F$ types. We know the gradation of  classical types has depth $\leq 2$. 
\subsection{Type $B$}
The gradation for $(B_\ell,\alpha_k)$ with $k\geq 2$ is given as follows, when $k=1$, it is Hermitian symmetric. For other $k$, we have the table
\begin{center}
	\captionof{table}{Roots for the gradation of $(B_\ell,\alpha_k)(k\geq 2)$}\label{Type B}
	\begin{tabular}{|p{9em}|p{9em}|p{9em}|}
		\hline
		$\Delta_0$ & $ \Delta_1$ & $\Delta_2$ \\  \hline
		$\pm(\lambda_i-\lambda_j), (i>k~\text{or}~j\leq k, 1\leq i<j\leq \ell);$ & $ \lambda_i-\lambda_j, (1\leq i\leq k <j\leq \ell);$ & $\lambda_i+\lambda_j, (1\leq i<j\leq k);$  \\
		$\pm\lambda_i, (k<i\leq \ell);$ & $\lambda_i, (1\leq i \leq k);$ &    \\
		$\pm(\lambda_i+\lambda_j), (k<i<j\leq \ell);$ & $\lambda_i+\lambda_j, (1\leq i\leq k<j\leq \ell);$ &\\
		\hline
	\end{tabular}
	
\end{center} 

When $k$ is even, we can obtain a set of maximal strongly orthogonal roots in $\Delta_2$ with long length, which is  $\{\lambda_1+\lambda_2,\lambda_3+\lambda_4,...,\lambda_{k-1}+\lambda_k\}$. Then we consider the elements in $\Delta_1$
\begin{enumerate}
	\item For $\lambda_i-\lambda_j, (1\leq i\leq k <j\leq \ell)$, we have a corresponding pair $(\lambda_i-\lambda_j,\lambda_{i+1}+\lambda_j)$ if $i$ is odd, $(\lambda_i-\lambda_j,\lambda_{i-1}+\lambda_j)$ if $i$ is even.
	\item For $\lambda_i$, we have a corresponding pair $(\lambda_i,\lambda_{i+1})$ if $i$ is odd, $(\lambda_i,\lambda_{i-1})$ if $i$ is even.
\end{enumerate}
Thus all roots in $\Delta_1$ have unique pairings.

\par When $k$ is odd, we still obtain a set of maximal strongly orthogonal roots in $\Delta_2$, which is  $\{\lambda_1+\lambda_2,\lambda_3+\lambda_4,...,\lambda_{k-2}+\lambda_{k-1}\}$. Then in $\Delta_1$, there is no pairing for $\{\lambda_k-\lambda_j, k<j\leq \ell\}$, $\{\lambda_k+\lambda_j, k<j\leq \ell\}$ and $\lambda_k$, for other roots there is a unique pairing. Then we are done.

\subsection{Type $C$}

For $(C_{\ell}, \alpha_k),(2\leq k\leq \ell-1)$ (for $k=1$ it is contact, for $k=\ell$ it is symmetric), we also list the roots as follows,

\begin{center}
	\captionof{table}{Roots for the gradation of $(C_\ell,\alpha_k)(2\leq k\leq \ell-1)$}\label{Type C}
	\begin{tabular}{|p{9em}|p{9em}|p{9em}|}
		\hline
		$\Delta_0$ & $ \Delta_1$ & $\Delta_2$ \\  \hline
		$\pm(\lambda_i-\lambda_j), (k<i~\text{or}~ j\leq k, 1\leq i<j\leq \ell);$ & $ \lambda_i-\lambda_j, (1\leq i\leq k <j\leq \ell);$ & $\lambda_i+\lambda_j, (1\leq i\leq j\leq k);$  \\
		$\pm(\lambda_i+\lambda_j), (k<i\leq j\leq \ell);$ & $\lambda_i+\lambda_j, (1\leq i\leq k<j\leq \ell);$ &    \\
		\hline
	\end{tabular}
\end{center}

\par We obtain a set of maximal strongly orthogonal roots in $\Delta_2$ with long length, which is $\{2\lambda_1, 2\lambda_2,2\lambda_3, 2\lambda_4,...,2\lambda_{k-1},2\lambda_k\}$. Then we consider the elements in $\Delta_1$
\begin{enumerate}
	\item For $\lambda_i-\lambda_j, (1\leq i\leq k <j\leq \ell)$, we have a corresponding pair $(\lambda_i-\lambda_j,\lambda_i+\lambda_j)$.
	\item For $\lambda_i+\lambda_j, (1\leq i\leq k <j\leq \ell)$, we have a corresponding pair $(\lambda_i+\lambda_j,\lambda_i-\lambda_j)$
\end{enumerate}
Then we are done since the pairing is unique if it exists.

\subsection{Type $D$}
\par For $(D_{\ell},\alpha_k)$ with $3\leq k\leq \ell-2$, we list the roots in each level. (If $k=1,\ell,\ell-1$, it is symmetric, if $k=2$, it is contact homogeneous.)

\begin{center}
	\captionof{table}{Roots for the gradation of $(D_\ell,\alpha_k)(3\leq k\leq \ell-2)$}\label{Type D}
	\begin{tabular}{|p{9em}|p{9em}|p{9em}|}
		\hline
		$\Delta_0$ & $ \Delta_1$ & $\Delta_2$ \\  \hline
		$\pm(\lambda_i-\lambda_j), (i-1\geq k~ \text{or} ~j\leq k, 1\leq i<j \leq \ell);$ & $ \lambda_i+\lambda_{\ell},(1\leq i\leq k);$ & $\lambda_i+\lambda_j,(1\leq i<j \leq  k);$  \\
		$\pm(\lambda_i+\lambda_j),(k+1\leq i<j \leq \ell-2);$ & $\lambda_i+\lambda_{\ell-1},(1\leq i\leq k);$ &    \\
		$\pm(\lambda_{\ell-1}+\lambda_{\ell});$& $\lambda_i+\lambda_j,(1\leq i\leq k<j \leq \ell-2);$&   \\
		$\pm(\lambda_i+\lambda_{\ell}), (k+1\leq i\leq \ell-2);$ & $\lambda_i-\lambda_j,(1\leq i\leq k<j \leq \ell);$& \\
		$\pm(\lambda_i+\lambda_{\ell-1}), (k+1\leq i\leq \ell-2);$& & \\
		\hline
	\end{tabular}
\end{center}

When $k$ is even, we obtain a set of maximal strongly orthogonal roots with long length, which is  $\{\lambda_1+\lambda_2,\lambda_3+\lambda_4,...,\lambda_{k-1}+\lambda_k\}$. We consider the elements in $\Delta_1$
\begin{enumerate}
	\item For $\lambda_i+\lambda_{\ell},(1\leq i\leq k)$, we have only one corresponding pair $(\lambda_i+\lambda_{\ell},\lambda_{i+1}-\lambda_{\ell})$ if $i$ is odd, $(\lambda_i+\lambda_{\ell},\lambda_{i-1}-\lambda_{\ell})$ if $i$ is even.
	\item For $\lambda_i+\lambda_{\ell-1},(1\leq i\leq k)$, we have only one corresponding pair $(\lambda_i+\lambda_{\ell-1},\lambda_{i+1}-\lambda_{\ell-1})$ if $i$ is odd, $(\lambda_i+\lambda_{\ell-1},\lambda_{i-1}-\lambda_{\ell-1})$ if $i$ is even.
	\item For $\lambda_i+\lambda_j, (1\leq i\leq k<j \leq \ell-2)$, we have only one corresponding pair $(\lambda_i+\lambda_j,\lambda_{i+1}-\lambda_j)$ if $i$ is odd, $(\lambda_i+\lambda_j,\lambda_{i-1}-\lambda_j)$ if $i$ is even.
	\item For $\lambda_i-\lambda_j,(1\leq i\leq k<j \leq \ell)$, we have only one corresponding pair $(\lambda_i-\lambda_j,\lambda_{i+1}+\lambda_j)$ if $i$ is odd, $(\lambda_i-\lambda_j,\lambda_{i-1}+\lambda_j)$ if $i$ is even.
\end{enumerate}

\par When $k$ is odd, we remain to have a set of maximal strongly orthogonal roots with long length, which is  $\{\lambda_1+\lambda_2,\lambda_3+\lambda_4,...,\lambda_{k-2}+\lambda_{k-1}\}$ in $\Delta_2$. In this case, there exists no pair for $\lambda_{k}+\lambda_{j}(k<j\leq \ell), \lambda_{k}-\lambda_{j}(k<j\leq \ell) $ in $\Delta_1$. But we can easily see that the pairing is unique if it exists.

\subsection{Type $E_6$}
\par For $(E_6, \alpha_k)$ with $k=3,5$ (Since when $k=1,6$ it is symmetric; when $k=2$ it is contact homogeneous; when $k=4$, the depth is three.)
we use the convention in \cite{MR0240238}, which starts from the root system of $E_8$. Let $\{e_i,1\leq i\leq 8\}$ be the standard basis of $\mathbb{R}^8$. Then let
\begin{align*}
&\alpha_1=\dfrac{1}{2}(e_8-e_7-e_6)+\dfrac{1}{2}(e_1-e_2-e_3-e_4-e_5), \alpha_2=e_1+e_2, \alpha_3=e_2-e_1,\\
&\alpha_4=e_3-e_2, \alpha_5=e_4-e_3, \alpha_6=e_5-e_4, \alpha_7=e_6-e_5, \alpha_8=e_7-e_6
\end{align*}
\par $E_6$ is generated by $\{\alpha_i,1\leq i\leq 6\}$. All  positive roots are $\{e_j\pm e_i, 1\leq i<j\leq 5\}\cup \{\dfrac{1}{2}(e_8-e_7-e_6+\sum_{i=1}^5(-1)^{v(i)}e_i), \sum_{i=1}^5v(i) \text{ is even}\}$
\par We list the roots as follows (when the associated root is $\alpha_3$):

\begin{center}
	\captionof{table}{Roots for the gradation of $(E_6,\alpha_3)$}\label{Type E_6}
	\begin{tabular}{|p{9em}|p{9em}|p{9em}|}
		\hline
		$\Delta_0$ & $ \Delta_1$ & $\Delta_2$ \\  \hline
		$\pm(e_j-e_i),(2\leq i<j\leq 5);$ & $ e_j-e_1,(1<j\leq 5);$ & $\dfrac{1}{2}(e_8-e_7-e_6+e_5+e_4+e_3+e_2+e_1);$  \\
		$\pm(e_j+e_1),(1<j\leq 5));$ & $e_j+e_i,(2\leq i<j\leq 5);$ &   $\dfrac{1}{2}(e_8-e_7-e_6+e_5+e_4+e_3-e_2-e_1)$; \\
		......&......& $\dfrac{1}{2}(e_8-e_7-e_6+e_5+e_4-e_3+e_2-e_1)$ \\ 
		&&$\dfrac{1}{2}(e_8-e_7-e_6+e_5-e_4+e_3+e_2-e_1)$ \\
		&&$\dfrac{1}{2}(e_8-e_7-e_6-e_5+e_4+e_3+e_2-e_1)$ \\
		\hline
	\end{tabular}
\end{center}

Roots in $\Delta_2$ can be reduced to only one root $\dfrac{1}{2}(e_8-e_7-e_6+e_5+e_4+e_3+e_2+e_1)=\alpha_1+\alpha_6+2\alpha_5+3\alpha_4+2\alpha_3+2\alpha_2$, then the unique pairing property is automatically satisfied. For $(E_6,\alpha_5)$, the process is the same since $\alpha_3$ and $\alpha_5$ are symmetric in the Dynkin diagram.

In some of the following cases, for simplicity we do not list the roots in $\Delta_0$ and only check that whether the difference of remaining roots in $\Delta_2$ is coincident with some difference of two roots in $\Delta_1$.

\subsection{Type $E_7$}
We consider $(E_7, \alpha_k),k=2,6$, since when $k=7$, it is symmetric; when $k=1$, it is contact homogeneous and when $k=2,6$, the depth is two.
\par $E_7$ is generated by $\{\alpha_i,1\leq i\leq 7\}$. All  positive roots are $\{e_j\pm e_i, 1\leq i<j\leq 6\}\cup\{e_8-e_7\}\cup \{\dfrac{1}{2}(e_8-e_7+\sum_{i=1}^6(-1)^{v(i)}e_i), \sum_{i=1}^6v(i) \text{ is odd}\}$.

We first consider $(E_7,\alpha_2)$. The coefficients for roots in $\Delta_2$ (written in terms of simple roots) are as follows, according to the order $ \alpha_1,\alpha_3,\alpha_4, \alpha_5, \alpha_6, \alpha_7, \alpha_2$.\[
\begin{gathered}
(2343212), (1343212), (1243212),
(1233212), (1232212), (1232112), (1232102)
\end{gathered}\]

We reduce the roots in $\Delta_2$ to only one root $ 2\alpha_1+3\alpha_3+4\alpha_4+3\alpha_5+2\alpha_6+\alpha_7+2\alpha_2=e_8-e_7$, then we finish.
\par For $(E_7,\alpha_6)$,The coefficients for roots in $\Delta_2$ (written in terms of simple roots) are as follows, according to the order $ \alpha_1,\alpha_3,\alpha_4, \alpha_5, \alpha_6, \alpha_7, \alpha_2$.
\[ \begin{gathered}
(2343212), (1343212), (1243212), (1233212),  (1232212)\\
(0122211), (1122211), (1222211), (1232211), (1233211)
\end{gathered}
\]
We obtain a set of maximal strongly orthogonal roots with long length in $\Delta_2$: $\{2\alpha_1+3\alpha_3+4\alpha_4+3\alpha_5+2\alpha_6+\alpha_7+2\alpha_2=e_8-e_7, \alpha_3+2\alpha_4+2\alpha_5+2\alpha_6+\alpha_7+\alpha_2=e_6+e_5\}$.  We cannot find two roots in $\Delta_1$ with difference $e_8-e_7-e_6-e_5$. Thus all the pairings must be unique if exist.

\subsection{Type $F_4$}
We consider $(F_4, \alpha_4)$ since for $(F_4, \alpha_1)$, it is contact homogeneous, for $(F_4,\alpha_2)$ and $(F_4,\alpha_3)$, the depth is at least three. 
\par Let $\{e_i,1\leq i\leq 4\}$ be the standard basis for $\mathbb{R}^4$. Then the positive simple roots can be written as $\alpha_1=e_2-e_3, \alpha_2=e_3-e_4, \alpha_3=e_4, \alpha_4=\dfrac{1}{2}(e_1-e_2-e_3-e_4)$. All positive roots are $\{e_i, 1\leq i\leq 4\} \cup\{e_i\pm e_j, 1\leq i<j\leq 4\} \cup \{\dfrac{1}{2}(e_1\pm e_2\pm e_3\pm e_4)\}$. For $(F_4, \alpha_4)$, the coefficients for roots in $\Delta_2$ (written in terms of simple roots) can be listed as follows, according to the order $\alpha_1, \alpha_2, \alpha_3, \alpha_4,$
\[\begin{gathered}
(0122), (1122), (1222), (1232), (1242), (1342), (2342)
\end{gathered}
\]
Then we have a set of maximal strongly orthogonal roots $\{2\alpha_1+3\alpha_2+4\alpha_3+2\alpha_4=e_1+e_2, \alpha_2+2\alpha_3+2\alpha_4=e_1-e_2\}$ in $\Delta_2$. We can not find two roots in $\Delta_1$ with difference $2e_2$. Thus all the pairings must be unique if exist.

\subsection{Type $E_8$}
We consider $(E_8,\alpha_1)$, since for $(E_8, \alpha_8)$ it is contact homogeneous and for other cases the depth is at least three.
\par  $E_8$ is generated by $\{\alpha_i,1\leq i\leq 8\}$. All positive roots are $\{e_j\pm e_i, 1\leq i<j\leq 8\}\cup \{\frac{1}{2}(e_8+\sum_{i=1}^7(-1)^{v(i)}e_i), \sum_{i=1}^7v(i) \text{ is even}\}$. 
\par For $(E_8,\alpha_1)$, the coefficients for roots in $\Delta_2$ (written in terms of simple roots) can be listed as follows, according to the order $\alpha_1, \alpha_3, \alpha_4, \alpha_5, \alpha_6, \alpha_7, \alpha_8, \alpha_2$
\[\begin{gathered}
(24654323), (24654313), (24654213), (24653213) \\ (24643213), (24543213),
(24543212), (23543213)\\
(23543212), (23443212),
(23433212), (23432212), (23432112), (23432102)
\end{gathered}
\]
Then we have a set of maximal strongly orthogonal roots  $\{2\alpha_1+4\alpha_3+6\alpha_4+5\alpha_5+4\alpha_6+3\alpha_7+2\alpha_8+3\alpha_2=e_8+e_7, 2\alpha_1+3\alpha_3+4\alpha_4+3\alpha_5+2\alpha_6+\alpha_7+2\alpha_2=e_8-e_7 \}$.  We cannot find two roots in $\Delta_1$ with difference $2e_7$  . Thus all the pairings must be unique if exist.

\section*{Acknowledgements}
The work was part of the author's PhD thesis. The author would like to thank his supervisor, Professor Ngaiming Mok, for his guidance.

	\bibliographystyle{alpha}      
\bibliography{research}   
\end{document}